\documentclass[12pt,reqno]{amsart}
\usepackage{cmap,mathtools}
\usepackage[T1]{fontenc}
\usepackage[utf8]{inputenc}
\usepackage[english]{babel}
\usepackage{xcolor}
\usepackage{amsfonts,amssymb,amsmath}
\usepackage{scrextend}
\usepackage{footnote,bigints}
\usepackage{array}
\usepackage{adjustbox}
\usepackage[labelformat=parens, labelfont={}]{subcaption}

\usepackage{float}
\usepackage{tikz}
\usepackage{amsthm}
\newtheorem{example}{Example}
\usepackage{caption}
\usepackage[margin=1in]{geometry}
\usepackage{graphicx}
\usepackage{pgfplots}
\pgfplotsset{compat=1.15}
\usepackage{mathrsfs}
\usetikzlibrary{arrows}
\definecolor{ududff}{rgb}{0.30196078431372547,0.30196078431372547,1.}
\usetikzlibrary{arrows}
 \usepackage{enumitem}
\allowdisplaybreaks
\usepackage[pagewise, mathlines]{lineno}

\usepackage{epstopdf}
\usepackage{a4wide}
\usepackage[colorlinks=true,linktocpage,pdfpagelabels,
bookmarksnumbered,bookmarksopen]{hyperref}
\definecolor{ForestGreen}{rgb}{0.1,0.6,0.05}
\definecolor{EgyptBlue}{rgb}{0.063,0.1,0.6}
\hypersetup{
	colorlinks=true,
	linkcolor=EgyptBlue,         
	citecolor=ForestGreen,
	urlcolor=olive
}

\usepackage[hyperpageref]{backref}
\newtheorem{theorem}{Theorem}[section]
\newtheorem{proposition}{Proposition}[section]
\newtheorem{definition}{Definition}[section]
\newtheorem{lemma}{Lemma}[section]
\newtheorem{remark}{Remark}[section]
\newtheorem{cor}{Corollary}[section]

\numberwithin{equation}{section}
\numberwithin{theorem}{section}
\numberwithin{equation}{section}
\numberwithin{theorem}{section}

\usepackage{ulem}
\usepackage{xcolor}
\usepackage[colorlinks=true,linktocpage,pdfpagelabels,
bookmarksnumbered,bookmarksopen]{hyperref}
\definecolor{ForestGreen}{rgb}{0.1,0.6,0.05}
\definecolor{EgyptBlue}{rgb}{0.063,0.1,0.6}
\hypersetup{
	colorlinks=true,
	linkcolor=EgyptBlue,         
	citecolor=ForestGreen,
	urlcolor=olive
}

\usepackage[hyperpageref]{backref}
\usepackage[foot]{amsaddr}
\usepackage{tikz}
\usetikzlibrary{graphs,arrows.meta} 
\DeclareUnicodeCharacter{2212}{-}
\title [Neumann eigenvalues on graphs]{Sharp bounds and monotonicity results for Neumann eigenvalues on graphs}
\author{Ashmita Singh\textsuperscript{1} \and Sheela Verma\textsuperscript{*}}
\address{$1$  Department of Mathematical Sciences, Indian Institute of Technology (BHU), Varanasi, India}
\email{ashmitasingh.rs.mat23@itbhu.ac.in}

\address{$*$  Corresponding author, Department of Mathematical Sciences, Indian Institute of Technology (BHU), Varanasi, India}
\email{sheela.mat@iitbhu.ac.in}
\DeclareUnicodeCharacter{2212}{-}
\usepackage{float}
\begin{document}
\begin{abstract}
In this article, we study sharp bounds for the Neumann eigenvalues of the Laplace operator on graphs. We first establish both lower and upper bounds for the second Neumann eigenvalue on simple graphs, and then derive monotonicity properties of Neumann eigenvalues on trees. In particular, we show that adding a vertex to a tree reduces the corresponding Neumann eigenvalues. We wish to emphasize that monotonicity results for Neumann eigenvalues on trees already exist in the literature, but our proof follows a fundamentally different approach. As a consequence of this monotonicity result, we provide an upper bound for the second Neumann eigenvalue and a lower bound for the largest Neumann eigenvalue on trees. Finally, we prove that under a diameter constraint on trees, the largest Neumann eigenvalue cannot be bounded from above. 
\end{abstract}
\maketitle
\section{Introduction}
\begin{figure}[H]
     \centering
     \resizebox{0.4\linewidth}{!}
{
\scalebox{0.4}
{
\begin{tikzpicture}[line cap=round,line join=round,>=triangle 45,x=1.0cm,y=1.0cm]

\draw (5*0.999048053014729,5*0.043619383008438) arc (2.5:357.5:5);
\draw (17-5*0.999048053014729,-5*0.043619383008438) arc (182.5:537.5:5);
\draw (5*0.999048053014729,5*0.043619383008438) -- (17-5*0.999048053014729,5*0.043619383008438);
\draw (5*0.999048053014729,-5*0.043619383008438) -- (17-5*0.999048053014729,-5*0.043619383008438);
\begin{scriptsize}

\draw[color=black] (12.5,0.2) node[font=\Huge] {$\epsilon$};
\end{scriptsize}
\end{tikzpicture}}
}
\captionof{figure}{\small Dumbbell domain.}
 \label{fig:dumbbell e.}
\end{figure}
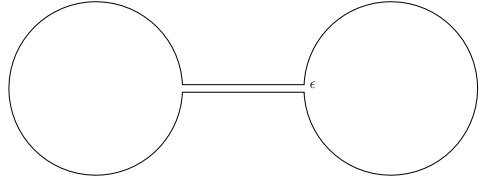

The Neumann boundary value problem was introduced in the nineteenth century by Carl Gottfried Neumann, a German mathematician, while studying potential theory. On a bounded smooth domain $M \subset \mathbb{R}^{n}$ with boundary $\partial M$, the Neumann eigenvalue problem is to find $(\lambda, f)$ satisfying 
\begin{align} \label{Neumann on mfd}
 \begin{cases}
 \Delta f(x) &= \lambda f(x), \quad  x \in M, \\
\frac{\partial f}{\partial n}(x) &= 0, \quad  x \in \partial M.
 \end{cases}   
\end{align}
 Here, the real number $\lambda$ is called a Neumann eigenvalue on $M$ and $f$ is called an eigenfunction corresponding to eigenvalue $\lambda$. The Neumann eigenvalues on $M$ are discrete and can be arranged as $0 = \lambda_1(M) <  \lambda_2(M) \leq  \lambda_3(M) \cdots \nearrow \infty$. The Neumann eigenvalues can also be interpreted as the frequencies of a freely vibrating membrane.
These eigenvalues are also interesting from spectral geometry point of view, as they contain information about the geometry of the domain under consideration. The problem of finding bounds for Neumann eigenvalues in terms of geometric quantities has been studied extensively. In general, it is not possible to compute Neumann eigenvalues explicitly except for special type of domains. However, some conclusions can be drawn about these eigenvalues in terms of the geometry of the domain. To see this, we consider the following examples:
\begin{example} \label{example rectangle}
 Consider a rectangular domain $M = [0, a] \times [0, b] \subset \mathbb{R}^2$. It is well known that the Neumann eigenvalues on the domain $M$ are given by
\begin{align*}
\lambda_{j,s} = \pi^2 \left( \frac{j ^2}{a^2} + \frac{s^2}{b^2} \right) \quad j, s = 0, 1, 2, \dots.
\end{align*}
If we arrange $\lambda_{j,s}$  in increasing order, then $\lambda_1(M) = \lambda_{0,0} = 0$ and value of each $\lambda_i(M), 2 \leq i < \infty$ depends on the values of $a$ and $b$.
\end{example}
\begin{example} \label{example dumbbell} 
Let $M_\epsilon$ denote a dumbbell domain in $\mathbb{R}^2$ (see Fig. \ref{fig:dumbbell e.}), consisting of two identical domains which are connected by a thin strip of width $\epsilon$. Then the second Neumann eigenvalue $\lambda_2(M_\epsilon)$ on $M_\epsilon$ approaches to $0$ as $\epsilon \to 0$. For more details, see \cite{colbois2006spectral}.
\end{example}

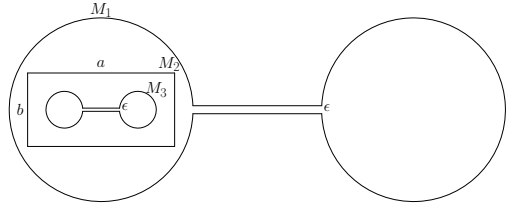
\begin{figure}[H]
     \centering
     \resizebox{0.45\linewidth}{!}
{
\scalebox{0.4}
{
\begin{tikzpicture}[line cap=round,line join=round,>=triangle 45,x=1.0cm,y=1.0cm]
\clip(-6.,-6.) rectangle (23.,6.);
\draw [line width=0.5pt] (-4.,2.)-- (4.,2.);
\draw [line width=0.5pt] (4.,2.)-- (4.,-2.);
\draw [line width=0.5pt] (4.,-2.)-- (-4.,-2.);
\draw [line width=0.5pt] (-4.,-2.)-- (-4.,2.);

\draw (-2+0.996194698092024,0.0871557427476582) arc (5:355:1);
\draw (2-0.996194698092024,-0.0871557427476582) arc (185:535:1);
\draw (-2+0.996194698092024,0.0871557427476582) -- (2-0.996194698092024,0.0871557427476582);
\draw (-2+0.996194698092024,-0.0871557427476582) -- (2-0.996194698092024,-0.0871557427476582);

\draw (5*0.999048053014729,5*0.043619383008438) arc (2.5:357.5:5);
\draw (17-5*0.999048053014729,-5*0.043619383008438) arc (182.5:537.5:5);
\draw (5*0.999048053014729,5*0.043619383008438) -- (17-5*0.999048053014729,5*0.043619383008438);
\draw (5*0.999048053014729,-5*0.043619383008438) -- (17-5*0.999048053014729,-5*0.043619383008438);
\begin{scriptsize}

\draw[color=black] (0,2.5) node[font=\Huge] {$a$};
\draw[color=black] (-4.4,0) node[font=\Huge] {$b$};
\draw[color=black] (0,5.45) node[font=\Huge] {$M_1$};
\draw[color=black] (3.7,2.5) node[font=\Huge] {$M_2$};
\draw[color=black] (3.0,1.1) node[font=\Huge] {$M_3$};
\draw[color=black] (12.3,0.1) node[font=\Huge] {$\epsilon$};
\draw[color=black] (1.3,0.16) node[font=\Huge] {$\epsilon$};
\end{scriptsize}
\end{tikzpicture}}
}
\captionof{figure}{\small Example to show that monotonicity property does not hold for Neumann eigenvalues.}
 \label{fig:dumbbell}
\end{figure}


From the above examples, we conclude the following well- known behaviors of Neumann eigenvalues on domains contained in the Euclidean space. 

\begin{enumerate}
    \item \textbf{Monotonicity Property:} In Fig. \ref{fig:dumbbell}, we consider domains $M_1, M_2, M_3 \subset \mathbb{R}^{2}$, where $M_1$ is a large dumbbell-shaped domain, $M_2$ is a rectangular domain, and $M_3$ is a small dumbbell-shaped domain, with $M_3 \subseteq M_2 \subseteq M_1$. For $M_1$ and $M_3$, by choosing the width of the connecting strip very small positive real number, $\lambda_2(M_1)$ and $\lambda_2(M_3)$ can be made arbitrarily small. Therefore, by keeping $M_2$ fixed and taking the strip width of $M_1$ and $M_3$ arbitrary small, using example \ref{example rectangle} and \ref{example dumbbell}, we get 
\begin{align*}
\lambda_2(M_1) \leq \lambda_2(M_2)\quad \text{and} \quad \lambda_2(M_3) \leq \lambda_2(M_2),
\end{align*}
where $\lambda_2$ represents the second Neumann eigenvalue. Hence, the monotonicity property does not hold for Neumann eigenvalues.
\item \textbf{Lower bound under diameter constraint:} From Example \ref{example dumbbell}, we observe that given any positive real numbers $D,c$, it is always possible to construct a dumbbell of diameter $D$ having strip width small enough such that its second Neumann eigenvalue is smaller than real number $c$. Therefore, it is not possible to establish a positive lower bound for the second Neumann eigenvalue on domains under a diameter constraint.
\end{enumerate}
  Our plan is to explore the above properties of Neumann eigenvalues on trees (defined in Section \ref{sec:prelim}) as trees are considered as a discrete analogue to Hadamard manifolds. A Hadamard manifold is defined as a simply connected Riemannian manifold with non positive sectional curvature.

Studying eigenvalue problems on graphs helps us understand spectrum of corresponding differential operators on manifolds. It is worth mentioning some classical articles where eigenvalues on graphs have been explored and that information has been utilized  to study the Laplacian spectrum on manifolds. Following the work of Buser \cite{buser1982note}, it has been recognized that discretizing a compact Riemannian manifold can be a highly effective approach for studying the spectrum of the Laplacian on that manifold. The idea of discretizing a manifold has also been used to study the Laplace and Steklov eigenvalues; see, for example, \cite{colbois2016steklov, mantuano2005discretization}. For further results related to this technique, see \cite{brooks2006combinatorial,brooks1986spectral, buser1984bipartition, chavel2001isoperimetric}.

Various studies have been done on the Laplace and Steklov eigenvalues on different graphs \cite{chung1997spectral,he2022upper,lin2024maximize,perrin2021isoperimetric,shi2022comparison,shi2022lichnerowicz}. However, Neumann eigenvalues have not been explored much. In \cite{shi2025comparisons}, the authors investigated Dirichlet, Neumann and Laplacian eigenvalues on weighted graphs. As applications, they derived Lichnerowicz-type, Fiedler-type and Friedman-type estimates for Dirichlet eigenvalues and Neumann eigenvalues. Monotonicity property of Neumann eigenvalues and Dirichlet eigenvalues have been obtained for some particular type of graphs in \cite{friedman1993some}. We acknowledge that the study of Steklov eigenvalues on graphs \cite{he2021steklov, perrin2019lower,yu2024monotonicity} has inspired parts of our approach, though this paper is exclusively concerned with Neumann eigenvalues.

In several situations, the Neumann eigenvalues on graphs are closely related with the unconstrained Laplace eigenvalues. Let $G$ be a simple and connected graph with boundary, and $G'$ be the graph obtained from $G$ by removing all boundary vertices. Let $\lambda_i(G)$ denote the $i$-th Neumann eigenvalue of $G$, and let $\mu_i(G')$ denote the $i$-th Laplace eigenvalue of $G'$ (with $i = 1, 2, \dots, |\Omega_G|$), where $\Omega_G$ denotes the set of interior vertices of $G$. Then the following properties hold
\begin{enumerate}
    \item Using the variational characterization of $\lambda_i(G)$ and $\mu_i(G')$, it can be shown that $\lambda_i(G) \geq \mu_i(G')$ for $i = 1, 2, \dots, |\Omega_G|$.
    \item If $G$ is a tree and $G'$ is obtained from $G$ by removing some boundary vertices (vertices having degree one), then $\lambda_i(G) = \mu_i(G')$. For details, see Remark \ref{rmk:Neumannproperty}.
\end{enumerate}

Throughout this article, $\lambda_i(G)$ denotes the $i^{th}$ Neumann eigenvalue of a graph with boundary, where $1 \leq i \leq k$ and $k$ is the number of interior vertices of $G$. 
We start by obtaining bounds for the second Neumann eigenvalue on graphs with boundary. More precisely, we present a sharp lower bound depending on the interior diameter and the number of interior vertices. 

\begin{theorem} \label{thm: lower bound}
Let $G$ be a connected graph with boundary, having interior diameter $d_\Omega$ and $k$ interior vertices. Then we have
\begin{align*}
\lambda_2 (G) \geq \frac{k}{(k-1)^2\, d_\Omega}.
\end{align*}
The bound is sharp up to a constant factor.
\end{theorem}
In addition, we establish an upper bound for the second Neumann eigenvalue on simple and connected graph $G$. Let $\Omega_G = \{v_1, v_2, \ldots, v_k\}$ be the set of interior vertices of $G$ which are numbered in such a way that $deg(v_1) = \displaystyle \min_{1 \leq i \leq k} deg(v_i)$ and $deg(v_2) = \displaystyle \min_{2 \leq i \leq k} deg(v_i)$. We denote $deg(v_1)=d_1$ and $deg(v_2)=d_2$. With these notations, we obtain the following theorem.
\begin{theorem}\label{the:upper bd on graph}
Let $G$ be a simple and connected graph with boundary. With the above notations, the second Neumann eigenvalue $\lambda_2(G)$ satisfies
\begin{align*}
\lambda_2(G) \leq
\begin{cases}
\frac{d_1+d_2}{2}, \quad \text{if } v_1 \sim v_2,\\
\frac{d_1+d_2}{2} + 1, \quad \text{if } v_1 \nsim v_2.
\end{cases}   
\end{align*}
\end{theorem}
We next study the behavior of Neumann eigenvalues on trees. Note that every tree is a simple graph, but the converse is not true. A tree is a connected graph that contains no cycles, whereas a simple graph may contain a cycle. For example, see Figures \ref{figPentagon},\ref{Fig:Pentagon 1} and \ref{Fig:Pentagon 2}; these are simple graphs but not trees. Crucially, The Neumann eigenvalues of a tree $G$ are the same as the Laplace eigenvalues of the tree $G'$, where $G'$ is obtained from $G$ by removing all its boundary vertices. The monotonicity property of Laplace eigenvalues was proved in \cite{MR1041245}, which in turn implies the monotonicity property of Neumann eigenvalues. However, in this article, we provide a different proof of the monotonicity result for Neumann eigenvalues on trees and establish the following result.
In particular, we show that Neumann eigenvalues satisfy an interlacing property when new vertices are added in such a way that the resulting graph is again a tree.
 
\begin{theorem} \label{thm:monotonicity}
Let $G_1$ and $G_2$ be two trees such that $G_1 \subset G_2$, that is, $G_2$ is obtained from $G_1$ by adding finitely many new edges. Let $k_1$ and $k_2$ denote the number of interior vertices in $G_1$ and $G_2$, respectively. Then for any $1 \leq i \leq k_1$, the $i^{\text{th}}$ Neumann eigenvalue satisfies $\lambda_i(G_2) \leq \lambda_i(G_1)  \leq \lambda_{i+k_2 - k_1}(G_2)$ . 
\end{theorem}

As a consequence of this monotonicity property, we obtain that on any tree the second Neumann eigenvalue is bounded above by $2$ and the largest Neumann eigenvalue is bounded below by $2$. 
\begin{cor} \label{cor:upper bound}
For any tree $T$ with $k$ interior vertices and diameter at least $3$, the second Neumann eigenvalue $\lambda_2$ and the largest Neumann eigenvalue $\lambda_k$ satisfy $$\lambda_2 (T) \leq 2 \leq \lambda_k(T).$$ 
\end{cor}
\begin{remark}\label{rem:mon on general graph}
Theorem \ref{thm:monotonicity} does not hold in general for arbitrary graphs. Consider graphs $G_1, G_2$ and $G_3$ as shown in Figures \ref{figPentagon}, \ref{Fig:Pentagon 1} and \ref{Fig:Pentagon 2}. Then Neumann eigenvalues on $G_1, G_2$ and $G_3$ are given as $\lambda_2(G_1)=1.5,\lambda_3(G_1)=2.5$; $\lambda_2(G_2)=1.5,\lambda_3(G_2)=3.5$ and $\lambda_2(G_3)=0.85380,\lambda_3(G_3)=1.6986$. Thus $G_1 \subset G_2$ and $G_1 \subset G_3$, however $\lambda_i(G_1) \leq \lambda_i(G_2)$ and $\lambda_i(G_1) \geq \lambda_i(G_3)$ for $i = 2,3$.
\end{remark}
\begin{figure}[h!]
  \centering
  \begin{subfigure}[t]{.3\textwidth}
    \centering
    \begin{tikzpicture}[scale=0.55, baseline=(current bounding box.south)]
      \clip(-2.5,-2.5) rectangle (2.5,2.5);
\draw [line width=1pt] (2,0)--(2,-2);
\draw [line width=1pt] (0,2)--(2,0);
\draw [line width=1pt] (0,2)--(-2,0);
\draw [line width=1pt] (-2,-2)--(2,-2);
\draw [line width=1pt] (-2,0)--(-2,-2);
\begin{scriptsize}
\draw [fill=black] (-2,-2) circle (1.5pt);
\draw (-2.2,-1.6) node {$v_5$};
\draw [fill=black] (2,-2) circle (3.5pt);
\draw (2.3,-1.6) node {$v_1$};
\draw [fill=black] (-2,0) circle (1.5pt);
\draw (-2.1,0.3) node {$v_4$};
\draw [fill=black] (2,0) circle (1.5pt);
\draw (2.2,0.3) node {$v_2$};
\draw [fill=black] (0,2) circle (3.5pt);
\draw (0.2,2.3) node {$v_3$};
\end{scriptsize}
    \end{tikzpicture}
    \caption{Graph $\mathbf{G_1}$ showing boundary (large dot) and interior (small dot) vertices.}
    \label{figPentagon}
  \end{subfigure}
  \hfill
  \begin{subfigure}[t]{.3\textwidth}
    \centering
    \begin{tikzpicture}[scale=0.55, baseline=(current bounding box.south)]
      \clip(-4.5,-2.5) rectangle (2.5,2.5);
\draw [line width=1pt] (2,0)--(2,-2);
\draw [line width=1pt] (0,2)--(2,0);
\draw [line width=1pt] (0,2)--(-2,0);
\draw [line width=1pt] (-2,-2)--(2,-2);
\draw [line width=1pt] (-2,0)--(-2,-2);
\draw [line width=1pt] (-4,-2)--(-2,-2);
\draw [line width=1pt] (-4,-2)--(-2,0);
\begin{scriptsize}
\draw [fill=black] (-2,-2) circle (1.5pt);
\draw (-1.8,-1.6) node {$v_5$};
\draw [fill=black] (2,-2) circle (3.5pt);
\draw (2.3,-1.6) node {$v_1$};
\draw [fill=black] (-2,0) circle (1.5pt);
\draw (-2.1,0.3) node {$v_4$};
\draw [fill=black] (2,0) circle (1.5pt);
\draw (2.2,0.3) node {$v_2$};
\draw [fill=black] (0,2) circle (3.5pt);
\draw (0.2,2.3) node {$v_3$};
\draw [fill=black] (-4,-2) circle (3.5pt);
\draw (-4.2,-1.6) node {$v_6$};
\end{scriptsize}
    \end{tikzpicture}
    \caption{Graph $\mathbf{G_2}$, obtained from $\mathbf{G_1}$ by adding a boundary vertex $v_6$.}
    \label{Fig:Pentagon 1}
  \end{subfigure}
  \hfill
  \begin{subfigure}[t]{.3\textwidth}
    \centering
    \begin{tikzpicture}[scale=0.55, baseline=(current bounding box.south)]
     \clip(-4.5,-2.5) rectangle (2.5,2.5);
\draw [line width=1pt] (2,0)--(2,-2);
\draw [line width=1pt] (0,2)--(2,0);
\draw [line width=1pt] (0,2)--(-2,0);
\draw [line width=1pt] (-2,-2)--(2,-2);
\draw [line width=1pt] (-2,0)--(-2,-2);
\draw [line width=1pt] (-4,-2)--(-2,-2);
\begin{scriptsize}
\draw [fill=black] (-2,-2) circle (1.5pt);
\draw (-1.8,-1.6) node {$v_5$};
\draw [fill=black] (2,-2) circle (3.5pt);
\draw (2.3,-1.6) node {$v_1$};
\draw [fill=black] (-2,0) circle (1.5pt);
\draw (-2.1,0.3) node {$v_4$};
\draw [fill=black] (2,0) circle (1.5pt);
\draw (2.2,0.3) node {$v_2$};
\draw [fill=black] (0,2) circle (3.5pt);
\draw (0.2,2.3) node {$v_3$};
\draw [fill=black] (-4,-2) circle (1.5pt);
\draw (-4.2,-1.6) node {$v_6$};
\end{scriptsize}
    \end{tikzpicture}
    \caption{Graph $\mathbf{G_2}$, obtained from $\mathbf{G_1}$ by adding an interior vertex $v_6$.}
    \label{Fig:Pentagon 2}
  \end{subfigure}
  \caption{Sequence of graphs $G_1$, $G_2$, and $G_3$ such that $G_1 \subset G_2$ and $G_1 \subset G_3$.}
  \end{figure}
Finally, we study the largest Neumann eigenvalue on the family of trees with fixed diameter. In particular, we prove the following result.

\begin{theorem} \label{thm:sup}
Let $\mathcal{G}$ be the collection of all trees with fixed diameter $D \geq 3$. Then the set $\{ \lambda_{max} (G) | G \in \mathcal{G} \}$ is not bounded above. Here, $\lambda_{max} (G)$ denotes the largest Neumann eigenvalue of $G$. 
\end{theorem}
 
This article is organised as follows. In Section \ref{sec:prelim}, we give some basic graph terminology to state the Neumann eigenvalue problem and present variational characterizations for Neumann eigenvalues. Bounds on the second Neumann eigenvalue on graphs have been obtained in Section \ref{Bounds for  the second Neumann eigenvalue}. In Section \ref{Some Auxiliary Results}, some auxiliary results are proved which are crucial for obtaining the monotonicity results. 
In Section \ref{Monotonicity property for Neumann eigenvalues}, we present the monotonicity results for Neumann eigenvalues on trees. Certain behaviors of Neumann eigenvalues on trees under diameter constraint have been discussed in Section \ref{Sharp Bound under Diameter constraint on trees}. Finally, in Section \ref{sec:concluding}, we state some open problems and potential extensions arising from this work.

\section{Preliminaries} \label{sec:prelim}
\subsection{Graph Terminology}
Let $G = (V, E)$ be a simple undirected graph with vertex set $V$ and edge set $E$. A graph is said to be simple if it is without loops and multiple edges. For any finite set $X$, we denote by $|X|$ its cardinality. In particular, $|V|$ denotes the number of vertices of $G$. Two vertices $x, y \in V$ are called adjacent, denoted by $x \sim y$ or $(x, y) \in E$, if there is an edge connecting $x$ and $y$. A path of length $\ell$ in $G$ is a sequence of distinct $(\ell+1)$ vertices $v_0,v_1, \cdots, v_\ell$ such that consecutive vertices are adjacent. The distance $d(x, y)$ between two vertices $x, y \in V$ is the length of the shortest path connecting them. The diameter $D$ of $G$ is the maximum distance among all pairs of vertices, namely, $D = \displaystyle \max_{x, y \in V} d(x, y)$. 
A graph $G$ is connected if for every pair of vertices $x, y \in V$, there exists a path connecting them. A closed path is called a cycle. A graph without cycle is called non-cylic. A tree is a non-cyclic connected graph. The degree of a vertex $v \in V$, denoted $\deg(v)$, is the number of edges incident to it, given by, $\deg(v) = |\{ u \in V : u \sim v \}|$. For any $S \subset V$,the vertex boundary of $S$, denoted $\delta S$, is the set of vertices in $S^c$ that are adjacent to $S$, given by $\delta S = \{ i \in S^c \mid i \sim j \text{ for some } j \in S \}$. A graph with boundary is a pair $(G,B)$ where $B \subset V(G)$ that satisfies $\delta (V \setminus B) = B$ and there is no edge connecting two vertices of $B$ i.e., $E(B, B) = \emptyset$. The set $B$ is called the boundary and $\Omega_{G} = V \setminus B$ is called the interior of the graph $G$. Note that $B=\delta\Omega_{G}$.

In this article, we consider Neumann eigenvalue problem on graph $G$ with total $n$ vertices and $k$ interior vertices, that is, $|V| = n$ and $|\Omega_{G}| = k$. If $G$ is a  tree, we consider $n \geq 3$ to exclude the case $\Omega_{G} = \emptyset$ or $k = 0$. All vertices of degree one in the tree $G$ are considered as boundary vertices of $G$. 

Now we state the Neumann eigenvalue problem on a graph.
\subsection{Neumann Eigenvalue Problem}
Let $G$ be a graph with interior $\Omega_{G}$ and boundary $\delta \Omega_{G}$. The vector space of real-valued functions on the set of vertices $V$ over $\mathbb{R}$, denoted by $\mathbb{R}^V$, is a finite-dimensional Euclidean space. For any $f, g \in \mathbb{R}^{V}$, we consider the inner product $\langle f, g \rangle_{\Omega_{G}} = \displaystyle \sum_{x \in \Omega_{G}} f(x)g(x)$. For $f \in \mathbb{R}^{V}$, the Laplace operator on the graph is defined as
\begin{align*}
(\Delta f)(x) = \sum_{\substack{y \in V \\ y \sim x}} (f(x) - f(y)),\quad x\in V.
\end{align*}

Within this framework, for any $f \in \mathbb{R}^V$ and $x \in \delta \Omega_{G}$, the outward normal derivative operator $\frac{\partial}{\partial \nu}: \mathbb{R}^V \to \mathbb{R}^{\delta\Omega_{G}}$, 
which maps $f$ to $\frac{\partial f}{\partial \nu}$, is expressed as
\begin{equation*}
\frac{\partial f}{\partial \nu}(x) = \sum_{\substack{
y \in{\Omega_{G}}\\ y \sim x}} (f(x) - f(y)) = (\Delta f)(x).
\end{equation*}

We now introduce the Neumann eigenvalue problem for the Laplace operator on graphs with boundary $(G,\delta \Omega_G)$. This problem is formulated as finding a nonzero function $ f \in \mathbb{R}^V$ and a real number $\lambda \in \mathbb{R}$ which satisfy the following system 
\begin{align} \label{Neumann}
 \begin{cases}
 \Delta f(x) &= \lambda f(x), \quad x \in \Omega_{G}, \\
\frac{\partial f}{\partial \nu}(x) &= 0, \quad  x \in \delta\Omega_{G}.
 \end{cases}   
\end{align}
The operator $\Lambda$ corresponding to the Neumann eigenvalue problem is given as follows. Define the map $\Lambda : \mathbb{R}^{\Omega_{G}} \to \mathbb{R}^{V}$ as
\begin{align*}
\Lambda(f)(x)=
\begin{cases}
f(x), & x \in \Omega_{G}, \\
\frac{1}{\operatorname{deg}(x)} \displaystyle\sum_{y \sim x} f(y), & x \in \delta\Omega_{G},
\end{cases}
\end{align*}
for any $f \in \mathbb{R}^{\Omega_{G}}$. 
It is straightforward to verify that for each $ x \in \delta\Omega_{G}$,
\begin{align*}
\frac{\partial (\Lambda(f))(x)}{\partial \nu} = 0.
\end{align*}
The Neumann eigenvalues on the pair $(G,\delta\Omega_{G})$ are same as the eigenvalues of the operator referred to as the Neumann Laplacian operator $\Delta^{N}: \mathbb{R}^{\Omega_{G}} \to \mathbb{R}^{\Omega_{G}}$, where
\begin{align*}
\Delta^{N} f = \Delta (\Lambda(f))|_{\Omega_{_G}}.
\end{align*}

Now, we show that $\Delta^{N}$ is a non-negative self-adjoint operator.
For any $f \in \mathbb{R}^{V}$, first, we define a map $d:\mathbb{R}^{V} \to A^{1}(G)$ as
\begin{align*}
df(x,y) := 
\begin{cases}
f(y) - f(x), & \{x,y\} \in E,\\
0, & \text{otherwise}.
\end{cases}
\end{align*}
Here $A^1(G)$ is the space of skew-symmetric functions $\alpha$ on $V \times V$ such that $\alpha(x, y) = 0$ when there is no edge connecting $x$ and $y$, equipped with the inner product $\langle \cdot , \cdot \rangle$
\begin{align}
\langle \alpha, \beta \rangle &= \frac{1}{2} \sum_{x, y \in V} \alpha(x, y)\beta(x, y) = \sum_{\{x, y\} \in E} \alpha(x, y)\beta(x, y). \label{eq:inner_product_A1}
\end{align}
Then for $f,g \in \mathbb{R}^V$, it follows easily by direct computation that
\begin{align*}
\langle \Delta^{N} f, g \rangle_{\Omega_{G}} =\langle d (f), d(g) \rangle.
\end{align*}
Thus the operator $\Delta^{N}$ is self-adjoint and non-negative; hence all its eigenvalues are real and non- negative. Therefore, for a graph $G$, there are $|\Omega_{G}|$ Neumann eigenvalues. These eigenvalues are non-negative and can be arranged in increasing order.
\begin{align*}
0 = \lambda_1 < \lambda_2 \leq \cdots \leq \lambda_{|\Omega_{G}|}.
\end{align*}
Here $\lambda_1 = 0$ because constant functions are eigenfunctions for the Neumann Laplacian operator, $\lambda_2 > 0$ due to the assumption that $G$ is connected.
In the next proposition, we present some well-known variational characterizations for $\lambda_i, 1 \leq i \leq |\Omega_G|$ and their proofs  are similar to the proof of the variational characterization of Laplace eigenvalues \cite{chung1997spectral}.
\begin{proposition}\label{prop:variational}
Let $G$ be a graph with $n$ vertices and $k$ interior vertices. For $i \leq k$, let $V_i = \operatorname{span}\{f^1, f^2, \ldots, f^i\}$, where $f^j \in \mathbb{R}^{n}$ is the eigenfunction corresponding to Neumann eigenvalue $\lambda_j, 1 \leq j \leq k$ on $G$. Then the $i$-th Neumann eigenvalue satisfies the following characterization
\begin{align}
\lambda_i(G) = R(f^i,G) 
&= \max_{\substack{f \in V_i \\ f \neq 0}} R(f,G) \label{R.Q.1}\\
&= \min_{\substack{f \in V_{i-1}^{\perp} \\ f \neq 0}} R(f,G)\label{R.Q.2}  \\
&= \min_{\substack{E \subset \mathbb{R}^n \\ \dim E = i}} \max_{\substack{f \in E \\ f \neq 0}} R(f,G)\label{R.Q.3}  \\
&= \max_{\substack{E \subset \mathbb{R}^n \\ \dim E = i-1}} \min_{\substack{f \in E^\perp \\ f \neq 0}} R(f,G)\label{R.Q.4}.
\end{align}
Here E is a subspace of $\mathbb{R}^n$ and $R(f,G) = \frac{\displaystyle\sum_{v_s \sim v_j } (f(v_s) - f(v_j))^2}{\displaystyle\sum_{v_s \in \Omega_G} f(v_s)^2}$ is the Rayleigh quotient. The sum $\displaystyle\sum_{v_s \sim v_j}$ is taken over all unordered pairs ${v_s,v_j}$ for which $v_s$ and $v_j$ are adjacent. The space $E^\perp$ is defined as
\begin{align*}
E^\perp  = \{ f \in \mathbb{R}^n \ | \ \langle f , g \rangle_{\Omega_G} = 0 \text{ for all } g \in E \}.
\end{align*}
\end{proposition}

\section{Bounds for the second Neumann eigenvalue} \label{Bounds for  the second Neumann eigenvalue}

We begin this section by discussing some properties of Neumann eigenfunctions on graphs. These properties will be used in calculating the sharp lower bound of second Neumann eigenvalue of graphs.

\begin{lemma} \label{lem: max min for eigenfunction}
Let $G$ be a connected graph with boundary, having $k$ interior vertices and $n$ total number of vertices. For $2 \leq i \leq k$, let $f^i = (f_1^i, f_2^i, \ldots, f_k^i,\ldots,f_n^i)$ be an eigenfunction corresponding to Neumann eigenvalue $\lambda_i(G)$. Then the function $f^i$ satisfies the following property.
\begin{enumerate}
    \item The maximum value and minimum value of $f^i_j$, $1 \leq j \leq n$, are attained at some interior vertices. \label{lem:max value at interior}
    \item The maximum value of $f^i_j$, $1 \leq j \leq n$ is strictly positive and the minimum value of $f^i_j$, $1 \leq j \leq n$ will be strictly negative. \label{lem:strict value at interior}
\end{enumerate}
\end{lemma}
\begin{proof}
$(1)$ If $f^i$ is an eigenfunction corresponding to the Neumann eigenvalue $\lambda_i(G)$. Then for any $x_j \in \delta\Omega_G$, the definition of Neumann problem \ref{Neumann} gives,
\begin{equation*}
\frac{\partial f^i}{\partial n}(x_j) = \sum_{\substack{
y \in{\Omega_{G}}\\ y \sim x_j}} (f^i(x_j) - f^i(y))=0.
\end{equation*}
Thus, $f^i(x_j)$ is equal to the average of function values at all the (interior) vertices adjacent to $x_j$. This ensures that the maximum and minimum values of $f^i_j$ are achieved at some interior vertices.\\

$(2)$ This follows easily from (\ref{lem:max value at interior}) and orthogonality condition $\displaystyle \sum_{v_j\in \Omega_{G}} f^i_j = 0$.
\end{proof}
We now introduce a family of graphs that will be used to demonstrate the sharpness of the lower bound of the second Neumann eigenvalue on graphs. 
\begin{definition}\label{def:G(s,a)}
Let $s \ge 1$ and $\alpha \ge 1$ be integers. Then the graph $G(s,\alpha)$ is constructed as follows. We start with a path $P_{2s+1} : v_0 \sim v_1 \sim \cdots \sim v_{2s+1}.$ For each vertex $v_j$ with $1 \le j \le 2s$, attach a complete graph $K_{\alpha-1}^{(j)}$ by connecting $v_j$ to every vertex of $K_{\alpha-1}^{(j)}$.
\end{definition}
\begin{figure}[ht]
\centering
\begin{tikzpicture}[
    interior/.style={circle, fill=black, inner sep=1pt},
    boundary/.style={circle, fill=black, inner sep=1.5pt},
    label/.style={font=\small}
]

\node[boundary] (v0) at (0,0) {};
\node[interior] (v1) at (2,0) {};
\node[interior] (v2) at (4,0) {};
\node           (dots) at (6,0) {\small$\cdots$};
\node[interior] (v2l2) at (8,0) {};
\node[interior] (v2l1) at (10,0) {};
\node[boundary] (v2l) at (12,0) {};

\node[label, above] at (v0) {$v_0$};
\node[label, above] at (v1) {$v_1$};
\node[label, above] at (v2) {$v_2$};
\node[label, above] at (v2l2) {$v_{2s-1}$};
\node[label, above] at (v2l1) {$v_{2s}$};
\node[label, above] at (v2l) {$v_{2s+1}$};

\draw (v0)--(v1)--(v2)--(dots)--(v2l2)--(v2l1)--(v2l);

\draw (2,-1.8) ellipse (0.9 and 0.5);
\draw (4,-1.8) ellipse (0.9 and 0.5);
\draw (8,-1.8) ellipse (0.9 and 0.5);
\draw (10,-1.8) ellipse (0.9 and 0.5);

\node[label] at (2,-1.8) {$K^{(1)}_{\alpha-1}$};
\node[label] at (4,-1.8) {$K^{(2)}_{\alpha-1}$};
\node[label] at (8,-1.8) {$K^{(2s-1)}_{\alpha-1}$};
\node[label] at (10,-1.8) {$K^{(2s)}_{\alpha-1}$};

\draw (v1)--(1.7,-1.3);
\draw (v1)--(2.3,-1.3);

\draw (v2)--(3.7,-1.3);
\draw (v2)--(4.3,-1.3);

\draw (v2l2)--(7.7,-1.3);
\draw (v2l2)--(8.3,-1.3);

\draw (v2l1)--(9.7,-1.3);
\draw (v2l1)--(10.3,-1.3);

\end{tikzpicture}
\caption{The graph $G(s,\alpha)$. Boundary vertices are shown as large dots, while interior vertices are shown as small dots.}
\label{fig:G(s,a)}
\end{figure}
Now we will give a lemma, the proof of which follows using a similar strategy as that of Lemma $2.4$ in \cite{lin2025estimates}. In our case, we choose a different test function in the variational characterization to get the desired bound.
\begin{lemma}\label{lem:pc-upper}
Let $G(s,\alpha)$ be a graph defined as above. Let $k$ be the number of interior vertices.
If $d_\Omega$ denotes the interior diameter of $G$, then
\begin{align}
\lambda_2(G) \leq \frac{12k}{(k-1)^2\,d_\Omega}.
\end{align}
Here, $k=2\alpha s$, and $d_\Omega=2s+1$ denotes the maximum distance among all pairs of interior vertices of $G$.
\end{lemma}
\begin{proof}
Let $v_0 \sim v_1 \sim \cdots \sim v_{2s+1}$ be the path of the graph $G=G(s,\alpha)$.
For $1 \le t \le 2s$, set
$F_t := \{v_t\} \cup V\!\left(K^{(t)}_{\alpha-1}\right).$
Then $k=2s\alpha$ and $d_\Omega=2s+1$.
Define a test function $f$ by
\begin{align} \label{eqn:test function sharpness}
f(v)=
\begin{cases}
s-(t-1), \quad &\text{if } v\in F_{t}, \qquad 1\le t\le s,\\
-(t-s), \quad &\text{if } v\in F_{t} , \qquad $s+1$\le t\le 2s,\\
s \quad &\text{if } v=v_0,\\
-s \quad &\text{if } v=v_{2s+1}.
\end{cases}
\end{align}
Then
\begin{align*}
\lambda_2(G) \leq R(f,G)
&= \frac{(s-1)+4+(s-1)}{{2\alpha \sum_{j=1}^s}j^2} \\
&= \frac{2s+2}{2\alpha (s(s+1)(2s+1)/6)} \\
&= \frac{6(s+1)}{\alpha s(s+1)(2s+1)} \\
&= \frac{6}{\alpha s\,d_\Omega}= \frac{12}{2\alpha\,s\,d_\Omega} \\
&=\frac{12}{k\,d_\Omega}=\frac{12k}{k^2\,d_\Omega}\leq \frac{12k}{(k-1)^2\,d_\Omega}.
\end{align*}
By \eqref{R.Q.2} of Proposition~\ref{prop:variational}, the eigenvalue
$\lambda_2(G)$ is characterized as the infimum of the Rayleigh quotient over all
functions $f$ orthogonal to the constant function. Note that the function $f$ defined in \eqref{eqn:test function sharpness} satisfies $\sum_{v_j \in \Omega_G} f(v_j)=0.$ Hence, $\lambda_2(G)\le R(f,G)$. Therefore, $\lambda_2(G) \leq \frac{12k}{(k-1)^2\,d_\Omega }$.
\end{proof}
Next, we provide a lower bound for the second Neumann eigenvalue on graphs. The idea of the proof is the same as used in Theorem $1$ of \cite{perrin2019lower}, (with slight modification for the Neumann condition) where a similar bound is derived for the second Steklov eigenvalue on graphs.\\
\begin{proof}[Proof of Theorem~\ref{thm: lower bound}]
Let $V = \{v_1, v_2, \ldots, v_n\}$ denote the vertices of $G$, and let $f = (f_1, f_2, \ldots, f_k,\ldots,f_n)$ be a normalized eigenfunction corresponding to eigenvalue $\lambda_2(G)$, where $f_j = f(v_j)$ for $1 \leq j \leq n$. Using Lemma \eqref{lem: max min for eigenfunction}, we obtain interior vertices $v_{k_1}$ and $v_{k_2}$ such that $\displaystyle \max_{v_j \in \Omega_G} f_j = f_{k_1}$ and $\displaystyle \min_{v_j \in \Omega_G} f_j = f_{k_2}$. Without loss of generality, we assume $f_{k_1} \geq |f_{k_2}|$.\\
From the normalization condition $\displaystyle \sum_{v_j \in \Omega_G} f_j^2 = 1$ we have
\begin{align*}
1 = \sum_{v_j\in \Omega_G} f_j^2 \leq k\, f_{k_1}^2 
\quad\implies\quad 
f_{k_1}^2 \geq \frac{1}{k} 
\quad\implies\quad 
f_{k_1} \geq \frac{1}{\sqrt{k}}.
\end{align*}
Using the condition $\displaystyle \sum_{v_j\in \Omega_G} f_j = 0$ we have
\begin{align*}
f_{k_1} + \sum_{\substack{v_j\in \Omega_G \\ j \neq {k_1}}} f_j = 0 
\quad\implies\quad 
-f_{k_1} = \sum_{\substack{v_j\in \Omega_G \\ j \neq {k_1}}} f_j.
\end{align*}
Since $f_{k_2}$ is the smallest value of $f_j$ for $1 \leq j \leq n$ on $\Omega_G$, we may write
\begin{align*}
-f_{k_1} = \sum_{\substack{v_j\in \Omega_G \\ j \neq {k_1}}} f_j \geq (k-1) f_{k_2}
\quad\implies\quad  f_{k_2} \leq -\frac{f_{k_1}}{k-1}.
\end{align*}
Substituting the earlier bound $f_{k_1} \geq \frac{1}{\sqrt{k}}$,
\begin{align*}
f_{k_2} \leq -\frac{f_{k_1}}{k-1} \leq -\frac{1}{(k-1)\sqrt{k}}.
\end{align*}
Therefore, we have
\begin{align*}
f_{k_1} - f_{k_2} \geq \frac{k}{(k-1) \sqrt{k}}.
\end{align*}
Since the graph is connected and the interior diameter of the graph is $d_\Omega$, there exists a path of length $\gamma \leq d_\Omega$ joining $v_{k_1}$ and $v_{k_2}$. We denote the $\gamma+1$ vertices of the path by $v_{i_1}, \dots, v_{i_{(\gamma+1)}}$, where $i_1 = {k_1}$ and $i_{(\gamma+1)} = {k_2}$. Using the Cauchy-Schwarz inequality, we obtain
\begin{align*}
\lambda_2 (G) = \sum_{v_s \sim v_j} (f_s - f_j)^2 \geq \sum_{t=1}^\gamma (f_{i_t} - f_{i_{t+1}})^2 \geq \frac{(f_{k_1} - f_{k_2})^2}{\gamma} \geq \frac{1}{d_\Omega} \left( \frac{k}{(k-1) \sqrt{k}} \right)^2 = \frac{k}{(k-1)^2\,d_\Omega}.
\end{align*}
The bound is sharp up to a constant factor, as follows from Lemma~\ref{lem:pc-upper}.
\end{proof}
\begin{remark}\label{rem:weighted}
Using similar technique, we can also establish a lower bound for the second Neumann eigenvalue on weighted graphs, which we plan to present in future work.
\end{remark}
\begin{proof}[Proof of Theorem~\ref{the:upper bd on graph}]
Let $\Omega_G$ denote the set of interior vertices of $G$.
The second Neumann eigenvalue $\lambda_2$ of the graph Laplacian is given by
\begin{align*}
\lambda_2(G)
= \min_{\substack{f \neq 0 \\
\sum_{v_s \in \Omega_G} f(v_s) = 0}}
\frac{\displaystyle\sum_{v_s \sim v_j } (f(v_s) - f(v_j))^2}{\displaystyle\sum_{v_s \in \Omega_{G}} (f(v_s))^2}.
\end{align*}

Define a test function $f$ by
\begin{align*}
f(v_i) =
\begin{cases}
1, & \text{if } v_i=v_1, \\
-1, & \text{if } v_i=v_2, \\
0, & \text{otherwise}.
\end{cases}
\end{align*}
The denominator of the Rayleigh quotient is
$$
\sum_{v_s \in \Omega_G} f(v_s)^2 = 1^2 + (-1)^2 = 2.
$$
We now compute the numerator.\\
\textbf{Case 1}: If $v_1 \sim v_2$.
$$
(f(v_1)-f(v_2))^2 = (1-(-1))^2 = 4.
$$
All remaining vertices adjacent to $v_1$ or $v_2$ contribute $1$ each.
Hence, the numerator equals
$$
(d_1-1) + (d_2-1) + 4 = d_1 + d_2 + 2.
$$
Therefore,
$$
\lambda_2 \leq \frac{d_1 + d_2 + 2}{2}
= \frac{d_1+d_2}{2} + 1.
$$
\textbf{Case 2: $v_1 \nsim v_2$.}\\
In this case, every vertex adjacent to $v_1$ or $v_2$ contributes $1$.
Thus, the numerator equals
$$
d_1 + d_2.
$$
Hence,
$$
\lambda_2 \leq \frac{d_1 + d_2}{2}.
$$
This completes the proof.
\end{proof}
\section{Some Auxiliary Results}\label{Some Auxiliary Results}
In this section, we provide some results that are used to establish the monotonicity property of Neumann eigenvalues on trees. Using the definition of linearly independent vectors, the following lemma follows easily.
\begin{lemma}\label{lem:Component equal}
Let $s \leq n$, consider a set of $s$ vectors $\{\mathbf{q}_1,\mathbf{q}_2,...,\mathbf{q}_s\}$ in $\mathbb{R}^{n+1}$ defined as
\begin{align*}
\mathbf{q}_1 &= (a_{1,1}, a_{1,2}, \ldots, a_{1,n}, a_{1,n+1}), \\
\mathbf{q}_2 &= (a_{2,1}, a_{2,2}, \ldots, a_{2,n}, a_{2,n+1}), \\
\vdots & \quad \vdots \\
\mathbf{q}_s &= (a_{s,1}, a_{s,2}, \ldots, a_{s,n}, a_{s,n+1}),
\end{align*}
such that $a_{t,n+1} = a_{t,j_0}$ for all $1\leq t \leq s$ and for some fixed $1 \leq j_0 \leq n$. Then the vectors $\{\mathbf{q}_1,\mathbf{q}_2,\ldots,\mathbf{q}_s\}$  are linearly independent if and only if the following set of vectors in $\mathbb{R}^n$ 
\begin{align*}
\mathbf{p}_1 &= (a_{1,1}, a_{1,2}, \ldots, a_{1,n}), \\
\mathbf{p}_2 &= (a_{2,1}, a_{2,2}, \ldots, a_{2,n}), \\
\vdots & \quad \vdots \\
\mathbf{p}_s &= (a_{s,1}, a_{s,2}, \ldots, a_{s,n}),
\end{align*}
obtained by omitting the $(n+1)$-th component of each $\mathbf{q}_t$, are linearly independent.
\end{lemma}
\begin{remark}\label{rmk:Neumannproperty}
Recall that in a tree $G$, all vertices of degree one are considered the boundary of the graph. Thus any boundary vertex is adjacent to exactly one interior vertex. Let $f$ be an eigenfunction corresponding to Neumann eigenvalue $\lambda(G)$ and $v_t$ be any boundary vertex adjacent to the interior vertex $v_s$. Then using the Neumann boundary condition on the boundary vertex $v_t$, we get $f(v_t) - f(v_s) = 0$ i.e., $f(v_t) = f(v_s).$ It can be easily concluded from this observation that if $G'$ is a graph obtained from a tree $G$ by removing all boundary vertices, then the Neumann eigenvalues of $G$ coincide with the Laplacian eigenvalues of $G'$. For a general graph, these two sets of eigenvalues are different.
\end{remark}

\begin{remark}\label{rmk:Neumannlaplace}
 In this article, we have also studied bounds on Neumann eigenvalues on simple graphs. Since we plan to study Neumann eigenvalue problems on general graphs in future work, this helps in building a basic framework for such studies. Due to this reason, our results on trees are stated in terms of Neumann eigenvalues in place of the Laplacian eigenvalues.
\end{remark}
Using Remark \ref{rmk:Neumannproperty}, a simplified expression of the variational characterization of Neumann eigenvalues on trees is established in the following proposition.
\begin{proposition} \label{prop:n to k}
Let $G$ be a tree on $n$ vertices, of which there are $k$ interior vertices. Then for $1 \leq i \leq k$, $\lambda_i(G)$ can also be characterized as 
\begin{align}\label{equivalent 1}
\lambda_i(G) = \min_{\substack{\tilde{E} \subset \mathbb{R}^k \\ \dim \tilde{E} = i}} \max_{\substack{\tilde{f} \in \tilde{E} \\ \tilde{f} \neq 0}} R({f},G).
\end{align}
Here, for any function $\tilde{f}$ defined on $\Omega_{G}$, its extension $f$ on $G$ is defined as follows
\begin{align}\label{f and tilde f}
    f(v_j) = 
    \begin{cases}
        \tilde{f}(v_j), \quad v_j \in \Omega_{G}, \\
        \tilde{f}(v_l), \quad v_j \in \delta\Omega_{G},
    \end{cases}
\end{align}
where $ v_l \in \Omega_G$ is such that $v_j \sim v_l.$
\end{proposition}
\begin{remark} \label{rmk: variartional chara}
Using Proposition \ref{prop:n to k},  we obtain the following simplified characterization of $\lambda_k(G),$
\begin{align*}
\lambda_k(G) =  \max_{\substack{\tilde{f} \in \mathbb{R}^k \\ \tilde{f} \neq 0}} R(f,G).
\end{align*}
Here $k$ is the number of interior vertices in $G$ and $f$ is the extension of $\tilde{f}$ on $G$ defined in \eqref{f and tilde f}.
\end{remark}

The following lemma is introduced to establish a one-to-one correspondence between the orthogonal complements of certain subspaces of $\mathbb{R}^{n}$ and $\mathbb{R}^{n+1}$, respectively. This will help us to use variational characterization \eqref{R.Q.4} in proof of Proposition \ref{prop:monotonicity2}.

\begin{lemma}
\label{one-one-correspondence}
Let \{$f^1, f^2, \ldots, f^{i-1}$\} be an orthogonal set of vectors in $\mathbb{R}^n$ with $f^\alpha = (f^\alpha_1, f^\alpha_2, \ldots, f^\alpha_n) \in \mathbb{R}^{n} \text{ for } 1 \leq\alpha \leq i-1$. Then $E =\operatorname{span}\{f^1, f^2, \ldots, f^{i-1}\}$ is an $(i-1)$-dimensional subspace of $\mathbb{R}^n$. Define subspace $\widetilde{E} \subset \mathbb{R}^{n+1}$ as
\begin{align*}
\widetilde{E} = \operatorname{span}\{\widetilde{f}^1, \widetilde{f}^2, \ldots, \widetilde{f}^{i-1}, u_l\},
\end{align*}
where $l$ is fixed natural number between  $1$ to $n$ and $\widetilde{f}^\alpha \in \mathbb{R}^{n+1}$ is defined by
\begin{align*}
\widetilde{f}^\alpha_j =
\begin{cases}
f^\alpha_j, & \text{for } j = 1, 2, \ldots, n \text{ with } j \neq l,\\
0, & \text{for } j = l,\\
f^\alpha_l, & \text{for } j = n+1,
\end{cases}
\end{align*}
and $u_l \in \mathbb{R}^{n+1}$ is defined as
$$
u_l = (0,0, \ldots, 0, \underbrace{1}_{l\text{-th}}, 0, \ldots, 0).
$$
Then the linear map $\phi: \mathbb{R}^{n} \to \mathbb{R}^{n+1}$ defined by
\begin{align*}
\phi(\mathbf{v}) = (v_1,v_2, \ldots, v_{l-1}, 0, v_{l+1}, \ldots, v_n, v_l) \text{ for } \mathbf{v} = (v_1,v_2, \ldots, v_n) \in \mathbb{R}^{n}
\end{align*}
is an isomorphism between $E^\perp$ and $\widetilde{E}^\perp$.
\end{lemma}
\begin{proof}
To establish a one-to-one correspondence between $E^\perp$ and $\widetilde{E}^\perp$, we construct a bijective linear map between them. Note that  $\dim(\widetilde{E}) = i$. Therefore, $\dim(\widetilde{E}^\perp) = (n+1) - i = (n - i + 1)$. Since $\dim(E^\perp) = (n - i + 1)$, the spaces $E^\perp$ and $\widetilde{E}^\perp$ have the same dimension and therefore they are isomorphic. 
Consider the linear map $\phi: \mathbb{R}^n \to \mathbb{R}^{n+1}$ defined by
\begin{align*}
\mathbf{\tilde{v}} = \phi(\mathbf{v}) = (v_1,v_2, \ldots, v_{l-1}, 0, v_{l+1}, \ldots, v_n, v_l),
\end{align*}
where $\mathbf{v} = (v_1, v_2, \ldots, v_n) \in \mathbb{R}^n$ and $\mathbf{\tilde{v}} = (\tilde{v}_1, \tilde{v}_2, \ldots, \tilde{v}_{n+1}) \in \mathbb{R}^{n+1}$. Note that $\widetilde{f}^\alpha = \phi ({f}^\alpha)$. Also, for any $\mathbf{v} \in \mathbb{R}^n$, $\mathbf{\tilde{v}}$ is orthogonal to $u_l$. The one-to-one correspondence between $E^\perp$ and $\widetilde{E}^\perp$ will be shown by restricting $\phi$ to $E^\perp$. For any $\mathbf{v}, \mathbf{w} \in \mathbb{R}^n$,
\begin{align} \label{ineq: orthogonality}
\langle \mathbf{\tilde{v}}, \mathbf{\tilde{w}}\rangle = \sum_{j=1}^{n+1} \tilde{v}_j \tilde{w}_j= \sum_{\substack{j=1}}^{n} \tilde{v}_j \tilde{w}_j + \tilde{v}_{n+1} \tilde{w}_{n+1} = \sum_{\substack{j=1 \\ j \neq l}}^n v_j w_j + v_l w_l = \sum_{j=1}^n v_j w_j = \langle \mathbf{v}, \mathbf{w} \rangle.
\end{align}
Thus, $\langle \mathbf{\tilde{v}}, \widetilde{f}^\alpha \rangle = \langle \mathbf{v}, f^\alpha \rangle$ and $\langle \mathbf{\tilde{v}}, u_l\rangle =0 $ for all $\mathbf{v} \in \mathbb{R}^{n}$. Hence, the linear map $\phi$ is an isomorphism between  $E^\perp$ and $\widetilde{E}^\perp$. This completes the proof. 
\end{proof}
\section{Monotonicity property for Neumann eigenvalues}\label{Monotonicity property for Neumann eigenvalues}
We are now prepared to present the monotonicity property of the Neumann eigenvalues on trees.
\begin{proposition}
\label{prop:monotonicity2}
Let $G_1$ be a finite tree with $n$ vertices and $k$ interior vertices, and let tree $G_2$ be obtained from $G_1$ by adding a new vertex adjacent to some boundary vertex of $G_1$. Then for $1 \leq i \leq k$,
\begin{enumerate}
\item $\lambda_i(G_1) \leq \lambda_{i+1}(G_2)$. \label{Ineq 2}
\item $\lambda_i(G_2) \leq \lambda_i(G_1)$. \label{Ineq 1}
\end{enumerate}
\end{proposition} 
\begin{proof}
Let $G_1$ be a finite tree having $n$ vertices and $k$ interior vertices with vertex set $V_1= \{ v_1,v_2,\ldots,v_n\}$. Let $G_2$ be obtained from $G_1$ by adding a vertex, say $v_{n+1}$, adjacent to some boundary vertex $v_l$ of $G_1$. The key observation is that the number of boundary vertices in tree $G_2$ is the same as the number of boundary vertices in $G_1$. However, $G_2$ has exactly one additional interior vertex than $G_1$ (See Figs.~\ref{fig:original graph} and \ref{fig:added to boundary}). 

$(1)$ We first show that $\lambda_i(G_1) \leq \lambda_{i+1}(G_2)$.
Let $f^1, f^2, \ldots, f^k \in \mathbb{R}^n$ be an orthogonal set of eigenfunctions of $G_1$ corresponding to the eigenvalues $\lambda_1(G_1) , \lambda_2(G_1) , \cdots,\lambda_k(G_1)$, respectively. Here,
\begin{align*}
f^\alpha = (f^\alpha_1, f^\alpha_2, \ldots, f^\alpha_n) = (f^\alpha(v_1), f^\alpha(v_2), \ldots, f^\alpha(v_n)) \in \mathbb{R}^{n} \text{ for all } 1 \leq \alpha \leq k.
\end{align*}
Define an ($i-1$)-dimensional subspace $E \subset \mathbb{R}^n$ as
\begin{align*}
E = \operatorname{span}\{f^1, f^2, \ldots, f^{i-1}\}.
\end{align*}
We now construct test functions in $\mathbb{R}^{n+1}$ corresponding to eigenvalue $\lambda_{i+1}(G_2)$. Let $v_l$ be the vertex where the new vertex $v_{n+1}$ is attached. For each vector $f^\alpha = (f^\alpha_1, f^\alpha_2, \ldots, f^\alpha_n) \in \mathbb{R}^n$ for $\alpha = 1, 2, \ldots, i-1$, define its extension $\widetilde{f}^\alpha \in \mathbb{R}^{n+1}$ as follows
\begin{align*}
\widetilde{f}^\alpha_j =
\begin{cases}
f^\alpha_j & \text{for } j = 1, 2, \ldots, n \text{ and } j \neq l,\\
0 & \text{for } j = l,\\
f^\alpha_l & \text{for } j = n+1.
\end{cases}
\end{align*}
Define the vector $u^l \in \mathbb{R}^{n+1}$ as the standard basis vector with a $1$ at the $l$-th position and $0$ elsewhere i.e.,
\begin{align*}
u^l = (0, \ldots, 0, \underbrace{1}_{l\text{-th}}, 0, \ldots, 0).
\end{align*}
Now, construct the $i$-dimensional subspace $\widetilde{E} \subset \mathbb{R}^{n+1}$ for $G_2$
\begin{align*}
\widetilde{E} = \operatorname{span}\{\widetilde{f}^1, \widetilde{f}^2, \ldots, \widetilde{f}^{i-1}, u^l\}.
\end{align*}
It follows from Lemma \ref{one-one-correspondence} that for any vector $g \in E^\perp$, its corresponding extension $\widetilde{g} \in \widetilde{E}^\perp$. By variational characterization of $\lambda_i(G_1)$ given in Proposition \ref{prop:variational}, equation \eqref{R.Q.4}
\begin{align}\label{unit vector_1}
    \lambda_i(G_1) = R(f^i, G_1) = \min_{\substack{f \in E^\perp \\ f \neq 0}} R(f,G_1). 
\end{align}
Note that for any $f \in E^{\perp}$,
\begin{align*}
R(f, G_1) = \frac{\displaystyle\sum_{v_s \sim v_j \in G_1} (f_s - f_j)^2}{\displaystyle\sum_{v_s \in \Omega_{G_1}} f_s^2}
\  \leq  \frac{\displaystyle\sum_{v_s \sim v_j \in G_1} (f_s - f_j)^2 + 2 f_l^2}{\displaystyle\sum_{v_s \in \Omega_{G_1}} f_s^2}  = \frac{\displaystyle\sum_{v_s \sim v_j \in G_2} (\widetilde{f}_s - \widetilde{f}_j)^2}{\displaystyle\sum_{v_s \in \Omega_{G_2}} \widetilde{f}_s^2}  = R(\widetilde{f}, G_2).
\end{align*}
Taking minimum over all nonzero $f \in E^{\perp}$, this inequality gives
\begin{align}\label{unit vector_2}
\lambda_i(G_1) = \min_{\substack{f \in E^\perp \\ f \neq 0}} R(f,G_1) \leq \min_{\substack{\widetilde{f} \in \widetilde{E}^\perp \\ \widetilde{f} \neq 0}} R(\widetilde{f}, G_2) \leq \lambda_{i+1}(G_2). 
\end{align}
Here the last inequality follows from the expression \eqref{R.Q.4} of variational characterization given in Proposition \ref{prop:variational}  of $\lambda_{i+1}(G_2)$. This completes the proof of the first part of the proposition.\\

$(2)$ Next, we prove the second part i.e., $\lambda_i(G_2) \leq \lambda_i(G_1)$ for all $1 \leq i \leq k$. Let $f^1, f^2, \ldots, f^i$ be orthogonal eigenfunctions corresponding to $\lambda_1(G_1), \lambda_2(G_1), \ldots, \lambda_i(G_1)$ respectively. Define $\widetilde{f}^\alpha$ on $G_2$ for $\alpha=1,2,\ldots,i$ by
\begin{align*}
\widetilde{f}^\alpha_j = \widetilde{f}^\alpha(v_j) = \begin{cases}
f^\alpha(v_j), & \text{for } j=1,2,\ldots,l\ldots, n, \\
f^\alpha(v_l), & \text{for } j=n+1.
\end{cases}
\end{align*}
Let $E = \operatorname{span}\{f^1,f^2,\ldots,f^i\}$ and $\widetilde{E} = \operatorname{span}\{\widetilde{f}^1,\widetilde{f}^2,\ldots,\widetilde{f}^i\}$. By Lemma \ref{lem:Component equal}, both $E$ and $\widetilde{E}$ are $i$-dimensional vector spaces. Then
\begin{align*}
\lambda_i(G_2) \leq \max_{\substack{\widetilde{f} \in \widetilde{E} \\ \widetilde{f} \neq 0}} R(\widetilde{f},G_2) \leq\max_{\substack{f \in E \\ f \neq 0}} R(f,G_1) = \lambda_i(G_1).
\end{align*}
The inequality $\lambda_i(G_2) \leq \max\limits_{\substack{\widetilde{f} \in \widetilde{E} \\ \widetilde{f} \neq 0}} R(\widetilde{f},G_2)$ holds by the variational characterization given in equation \eqref{R.Q.3} and the fact that $\widetilde{E}$ is an $i$-dimensional subspace of $\mathbb{R}^{n+1}$. The inequality
 $\max\limits_{\substack{\widetilde{f} \in \widetilde{E} \\ \widetilde{f} \neq 0}} R(\widetilde{f},G_2) \leq \max\limits_{\substack{f \in E \\ f \neq 0}} R(f,G_1)$ holds, since the numerator remains unchanged in both cases, while the denominator of $R(\widetilde{f},G_2)$ is larger than the denominator of $R(f,G_1)$, since $G_2$ has one extra interior vertex than $G_1$. Finally $\max\limits_{\substack{f \in E \\ f \neq 0}} R(f,G_1) = \lambda_i(G_1)$ is true by Proposition \ref{prop:variational}. Therefore, we conclude that
\begin{align*}
\lambda_i(G_2) \leq \lambda_i(G_1).
\end{align*}
\end{proof}
\begin{proposition}
\label{prop:monotonicity1}
Let $G_1$ be a finite tree with total $n$ vertices and $k$ interior vertices, and tree $G_2$ be obtained from $G_1$ by adding a new vertex adjacent to some interior vertex of $G_1$. Then for $1 \leq i \leq  k$,
\begin{align}
\lambda_i(G_1) = \lambda_i(G_2). 
\end{align}
\end{proposition} 
\begin{proof}
Proof follows easily by using Remark \ref{rmk:Neumannproperty}.
\end{proof}

\vspace{-1em} 
\begin{figure}[h]
\centering
\subfloat[Graph \textbf{G}$_\mathbf{1}$   showing boundary vertices (large dots) and interior vertices (small dots).]{\label{fig:original graph}\scalebox{1}{\begin{tikzpicture}[line cap=round,line join=round,>=stealth,x=0.9cm,y=0.9cm]
        \clip(-1.2,-0.5) rectangle (6.2,3.2);
        \draw [line width=1pt] (1.,1.)-- (2.,1.);
        \draw [line width=1pt] (2.,1.)-- (3.,1.);
        \draw [line width=1pt,] (3.,1.)-- (4.,1.);
        \draw [line width=1pt,] (4.,1.)-- (5.,1.);
        \draw [line width=1pt,] (2.,1.)-- (2.,2.);
        \draw [line width=1pt,] (4.,1.)-- (4.,2.);
        \draw [line width=1pt,] (4.,1.)-- (4.6,0.3);
        \begin{scriptsize}
        \draw [fill=black] (1.,1.) circle (3pt);
        \draw[color=black] (1.,0.7) node {$v_1$};
        \draw [fill=black] (2.,1.) circle (1.5pt);
        \draw[color=black] (2.,0.7) node {$v_2$};
        \draw [fill=black] (3.,1.) circle (1.5pt);
        \draw[color=black] (3.,0.7) node {$v_3$};
        \draw [fill=black] (4.,1.) circle (1.5pt);
        \draw[color=black] (4.,0.7) node {$v_4$};
        \draw [fill=black] (5.,1.) circle (3pt); 
        \draw[color=black] (5.,0.7) node {$v_5$};
        \draw [fill=black] (4.6,0.3) circle (3pt);
         \draw[color=black] (4.6,0.05) node {$v_6$};
        \draw [fill=black] (2.,2.) circle (3pt); 
         \draw[color=black] (2.4,2) node {$v_7$}; 
         \draw [fill=black] (4.,2.) circle (3pt); 
         \draw[color=black] (4.4,2) node {$v_8$};
        \end{scriptsize}
    \end{tikzpicture}}}
\hfill
\subfloat[Graph \textbf{G}$_\mathbf{2}$, obtained from \textbf{G}$_\mathbf{1}$ by adding vertex $v_9$ to boundary vertex $v_5$.]{ \label{fig:added to boundary}\scalebox{1}{\begin{tikzpicture}[line cap=round,line join=round,>=stealth,x=0.9cm,y=0.9cm]
        \clip(-1.2,-0.5) rectangle (8.2,3.2);
        \draw [line width=1pt] (1.,1.)-- (2.,1.);
        \draw [line width=1pt] (2.,1.)-- (3.,1.);
        \draw [line width=1pt,] (3.,1.)-- (4.,1.);
        \draw [line width=1pt,] (4.,1.)-- (5.,1.);
        \draw [line width=1pt,] (2.,1.)-- (2.,2.);
        \draw [line width=1pt,] (4.,1.)-- (4.,2.);
        \draw [line width=1pt,] (4.,1.)-- (4.6,0.3);
        \draw [line width=1pt,dashed] (5.,1.)-- (6.,1.);
        \begin{scriptsize}
        \draw [fill=black] (1.,1.) circle (3pt);
        \draw[color=black] (1.,0.7) node {$v_1$};
        \draw [fill=black] (2.,1.) circle (1.5pt); 
        \draw[color=black] (2.,0.7) node {$v_2$};
        \draw [fill=black] (3.,1.) circle (1.5pt);
        \draw[color=black] (3.,0.7) node {$v_3$};
        \draw [fill=black] (4.,1.) circle (1.5pt);
        \draw[color=black] (4.,0.7) node {$v_4$};
        \draw [fill=black] (5.,1.) circle (1.5pt); 
        \draw[color=black] (5.,0.7) node {$v_5$};
        \draw [fill=black] (4.6,0.3) circle (3pt); 
         \draw[color=black] (4.6,0.05) node {$v_6$};
        \draw [fill=black] (2.,2.) circle (3pt); 
         \draw[color=black] (2.4,2) node {$v_7$}; 
         \draw [fill=black] (4.,2.) circle (3pt); 
         \draw[color=black] (4.4,2) node {$v_8$};
         \draw [fill=black] (6.,1.) circle (3pt);
         \draw[color=black] (6.,0.7) node {$v_9$};

        \end{scriptsize}
        \end{tikzpicture}}}
\caption{Examples.}
\label{fig}
\end{figure}
Using Propositions \ref{prop:monotonicity2} and \ref{prop:monotonicity1}, we conclude the following monotonicity property of the Neumann eigenvalues on trees.\\
\begin{proof}[Proof of Theorem \ref{thm:monotonicity}]
The first inequality, $\lambda_i(G_2) \leq \lambda_i(G_1)$ 
follows directly from Proposition \ref{prop:monotonicity1} and part \eqref{Ineq 1} of Proposition \ref{prop:monotonicity2}. For the inequality $\lambda_i(G_1)  \leq \lambda_{i+k_2 - k_1}(G_2)$, we use Proposition \ref{prop:monotonicity1} and part \eqref{Ineq 2} of Proposition \ref{prop:monotonicity2}. Let $G_2$ be obtained from $G_1$ in $j$ steps, where in each step we add exactly one edge, i.e., $G_1 = G_{1,1} \subset G_{1,2} \subset G_{1,3} \subset \cdots \subset G_{1,j+1} = G_{2}$, where $G_{1, s+1}$ has exactly one extra edge and vertex (interior or boundary) from $G_{1, s}$ for each $1 \leq s \leq j$. Fix some $i$ with $1 \leq i \leq k_1$, then for each $1 \leq s \leq j$,
\begin{align*}
\lambda_i(G_{1,s}) \leq
\begin{cases}
\lambda_i(G_{1,s+1}) , &\quad \text{ If } G_{1,s+1} \text { has one extra boundary vertex,} \\
\lambda_{i+1}(G_{1,s+1}) , &\quad \text{ If } G_{1,s+1} \text { has one extra interior vertex.} 
\end{cases}
\end{align*}
Note that in the process of obtaining $G_2$ from $G_1$ by adding one vertex at each step, there are exactly $(k_2 - k_1)$ steps in which an interior vertex is added. Combining all the above facts, we conclude that $\lambda_i(G_1)\leq\lambda_{i+k_2 - k_1}(G_2)$.
\end{proof}
As a consequence of the above theorem, the following result is obtained, providing an upper bound for $\lambda_2$ and a lower bound for $\lambda_k$ for trees having $k$ interior vertices. \\
\begin{proof}[Proof of Corollary~\ref{cor:upper bound}]
Note that any tree $T$ with diameter at least $3$ contains path $P_3$ of diameter three. It is easy to compute that the path graph of diameter 3 has two Neumann eigenvalues, namely $\lambda_1(P_3) = 0$ and $\lambda_2(P_3) = 2$. Then $\lambda_2(T) \leq \lambda_2(P_3) \leq \lambda_{2+k-2}(T) = \lambda_{k}(T)$. Since $\lambda_2(P_3) = 2,$ we get the desired result.
\end{proof}
\section{Sharp Bound under Diameter Constraint on trees}\label{Sharp Bound under Diameter constraint on trees}
Now, we provide a theorem demonstrating that, under a diameter constraint, the set of largest Neumann eigenvalues on the family of trees is unbounded.
\begin{proof}[Proof of Theorem \ref{thm:sup}]
\begin{figure}[h]
    \centering
    \resizebox{0.75\linewidth}{!}
    {
        \begin{tikzpicture}[line cap=round,line join=round,>=triangle 45,x=1.0cm,y=1.0cm]
          \draw[line width=1.0pt] (0.,0.) -- (1.1490909090909092,0.);
            \draw[line width=1.0pt] (1.1490909090909092,0.) -- (2.38,0.);
             \draw[line width=1.0pt] (0.,0.) -- (0.62,0.76);
            \draw[line width=1.0pt] (0.62,0.76) -- (1.09,1.42);  
             \draw[line width=1.0pt] (-1.14,0.) -- (-2.42,-0.02);
            \draw[line width=1.0pt] (0.,0.) -- (-1.14,0.);
            \draw[line width=1.0pt] (0.,0.) -- (-0.5,0.74);
            \draw[line width=1.0pt] (-0.5,0.74) -- (-1.02,1.4);

            \draw[line width=1.0pt] (0.,0.) -- (-0.54,-0.68);
            \draw[line width=1.0pt] (-0.54,-0.68) -- (-1.08,-1.34);

            \draw[line width=1.0pt] (0.,0.) -- (0.56,-0.66);
            \draw[line width=1.0pt] (0.56,-0.66) -- (1.18,-1.32);

            \draw[line width=1.0pt] (2.38,0.) -- (3.5854545454545463,0.);
            \draw[line width=1.0pt] (3.5854545454545463,0.) -- (4.84,0.);
            \draw[line width=1.0pt,dashed] (4.84,0.) -- (6.,0.);
            \draw[line width=1.0pt] (6.,0.) -- (7.185454545454551,-0.010909090909098557);
            \draw[line width=1.0pt,dashed] (-0.5,0.74) -- (-1.14,0.);
            \draw[line width=1.0pt,dashed] (-1.02,1.4) -- (-2.42,-0.02);
            \draw[line width=1.0pt,dashed] (-1.14,0.) -- (-0.54,-0.68);
            \draw[line width=1.0pt,dashed] (-2.42,-0.02) -- (-1.08,-1.34);
             \begin{scriptsize}
                \draw[fill=black] (0.,0.) circle (1.5pt);
                \draw (0.30,0.10) node[font=\scriptsize] {$o$};
               \draw[fill=black] (1.1490909090909092,0.) circle (1.5pt);
                \draw (1.276363636363639,0.34363636363635597) node[font=\scriptsize] {$c_1$};
                \draw[fill=black] (2.38,0.) circle (1.5pt);
                \draw (2.512727272727276,0.34363636363635597) node[font=\scriptsize] {$c_2$};

                \draw[fill=black] (0.62,0.76) circle (1.5pt);
                \draw (0.7490909090909114,1.0890909090909013) node[font=\scriptsize] {$b_1$};

               \draw[fill=black] (1.09,1.42) circle (3pt); 
                \draw (1.1490909090909116,1.7618181818181742) node[font=\scriptsize] {$a_1$};
               
                \draw[fill=black] (-1.14,0.) circle (1.5pt);
                \draw (-1.014545454545453,0.34363636363635597) node[font=\scriptsize] {$b_l$};
                \draw[fill=black] (-2.42,-0.02) circle (3pt); 
                \draw (-2.2872727272727267,0.3254545454545378) node[font=\scriptsize] {$a_l$};
                \draw[fill=black] (-0.5,0.74) circle (1.5pt);
                \draw (-0.3781818181818165,1.0709090909090833) node[font=\scriptsize] {$b_2$};
                \draw[fill=black] (-1.02,1.4) circle (3pt); 
                \draw (-0.8872727272727259,1.7436363636363559) node[font=\scriptsize] {$a_2$};

                \draw[fill=black] (-0.54,-0.68) circle (1.5pt);
                \draw (-0.4145454545454529,-0.34727272727273495) node[font=\scriptsize] {$b_j$};
                \draw[fill=black] (-1.08,-1.34) circle (3pt);
                \draw (-0.96,-1.0018181818181895) node[font=\scriptsize] {$a_j$};

                \draw[fill=black] (0.56,-0.66) circle (1.5pt);
                \draw (0.6945454545454568,-0.3290909090909167) node[font=\scriptsize] {$b_{j+1}$};
                \draw[fill=black] (1.18,-1.32) circle (3pt);
                \draw (1.3127272727272754,-0.9836363636363712) node[font=\scriptsize] {$a_{j+1}$};

                \draw[fill=black] (3.5854545454545463,0.) circle (1.5pt);
                \draw (3.7127272727272764,0.34363636363635597) node[font=\scriptsize] {$c_3$};
                \draw[fill=black] (4.84,0.) circle (1.5pt);
                \draw (4.967272727272732,0.34363636363635597) node[font=\scriptsize] {$c_4$};
                 \draw[fill=black] (6.,0.) circle (1.5pt); 
                \draw (6.130909090909096,0.34363636363635597) node[font=\scriptsize] {$c_{D-3}$};
                \draw[fill=black] (7.185454545454551,-0.010909090909098557) circle (3pt); 
                \draw (7.35454545454551,0.34363636363635597) node[font=\scriptsize] {$c_{D-2}$};
            \end{scriptsize}
        \end{tikzpicture}
   } 
    \caption{Tree $G_{j+1}$ with diameter $D$, where $o, b_1, b_2, \dots, b_{j+1}, c_1, c_2, \dots, c_{D-3}$ are interior vertices, and $a_1, a_2\dots, a_{j+1}, c_{D-2}$ are leaves.
}
    \label{fig:tree_diameter_D}
\end{figure}
To prove this theorem, we construct a sequence of trees $G_{j+1}, j \in \mathbb{N}$ (see Fig.\ref{fig:tree_diameter_D}) as follows:
Start with a path $P$ on $D-1$ vertices, labeled as
$o, c_1, c_2, \cdots, c_{D-2}$. To the vertex $o$, attach $j+1$ paths of length $2$ each. That is, for each $t = 1, 2, \ldots, j+1$, we introduce two vertices $a_t$ and $b_t$ such that $o \sim b_t$  \text{ and } $b_t \sim a_t.$
The resulting tree $G_{j+1}$ has diameter $D$ with $j+D-1$ interior vertices and $j+2$ boundary vertices. Now we show that the maximum Neumann eigenvalue of $G_{j+1}$ converges to infinity as $j$ approaches infinity. Using Remark \ref{rmk: variartional chara}, we write
\begin{align*}
\lambda_{\max}(G) = \max_{\substack{\tilde{f} \in \mathbb{R}^{\Omega_{G}} \\ \tilde{f} \neq 0}}  R(f,G) =  \frac{\displaystyle\sum_{v_s \sim v_j } (f(v_s) - f(v_j))^2}{\displaystyle\sum_{v_s \in \Omega_{G}} f(v_s)^2}.
\end{align*}
Define the test function $f$ on $G_{j+1}$ as
\begin{align*}
f(v) = \begin{cases} 
-(j + D - 2), & \text{if } v = o,\\
1, & \text{otherwise }.
\end{cases}
\end{align*}

Then
\begin{align*}
\lambda_{\max} (G_{j+1}) \geq \frac{(1+(j+D-2))^2(j+2)}{1^2 \cdot ({j+1})+1^2 \cdot {(D-3)} + (j+D-2)^2 \cdot 1}.
\end{align*}
Simplifying this expression, we get
\begin{align*}
\lambda_{\max} (G_{j+1}) &\geq \frac{(j+D-1)^2(j+2)}{(j+D-2) + (-(j+D-2))^2} \\
&= \frac{(j+2)(j+D-1)^2}{(j+D-2)(1 + j + D - 2)} \\
&= \frac{(j+2)(j+D-1)}{j+D-2} \\
&= \frac{j^2(1+\frac{2}{j})(1+\frac{D-1}{j})}{j(1+\frac{D-2}{j})} \\
&\to \infty \quad \text{as } j \to \infty.
\end{align*}
Hence, the proof is complete.
\end{proof}

\section{Concluding Remarks}\label{sec:concluding}
We conclude this article by listing some open problems and potential extensions arising from this work.

\begin{enumerate}
    \item Theorem \ref{thm:monotonicity} establishes an interlacing property for Neumann eigenvalues when a tree is enlarged by adding vertices. However, as shown in Remark \ref{rem:mon on general graph}, this monotonicity property does not hold for arbitrary graphs. A natural question is: \textit{To characterize the largest class of graphs (containing trees) for which the Neumann eigenvalues satisfy an interlacing property similar to Theorem \ref{thm:monotonicity}.} 

    \item The variational characterization technique to prove monotonicity result can potentially be extended to prove monotonicity results for other spectral problems on graphs.

    \item A systematic study of Neumann eigenvalues on weighted graphs with boundary remains to be explored further.
\end{enumerate}

\textbf{Acknowledgements} A. Singh is thankful to IIT (BHU) Varanasi for providing financial support in the form of a scholarship. The corresponding author S. Verma acknowledges the project grant provided by CSIR-ASPIRE sanction order no. 25WS(011)/2023-24/EMR-II/ \\ ASPIRE and SERB-SRG sanction order No. SRG/2022/002196.

\textbf{Data Availability} No data, models, or code were generated or used during the study.

\section*{Declarations}
\textbf{Conflicts of Interest} The author declared no potential conflicts of interests with respect to this article.

\bibliographystyle{plain}
\bibliography{references}

@book {chung1997spectral,
    AUTHOR = {Chung, Fan R. K.},
     TITLE = {Spectral graph theory},
    SERIES = {CBMS Regional Conference Series in Mathematics},
    VOLUME = {92},
 PUBLISHER = {Conference Board of the Mathematical Sciences, Washington, DC; by the American Mathematical Society, Providence, RI},
      YEAR = {1997},
     PAGES = {xii+207},
      ISBN = {0-8218-0315-8},
   MRCLASS = {58G99 (05C50 35P05 46N20 47N20)},
  MRNUMBER = {1421568},
MRREVIEWER = {Robert\ Brooks},
}

@article {buser1982note,
    AUTHOR = {Buser, Peter},
     TITLE = {A note on the isoperimetric constant},
   JOURNAL = {Ann. Sci. \'Ecole Norm. Sup. (4)},
  FJOURNAL = {Annales Scientifiques de l'\'Ecole Normale Sup\'erieure. Quatri\`eme S\'erie},
    VOLUME = {15},
      YEAR = {1982},
    NUMBER = {2},
     PAGES = {213--230},
      ISSN = {0012-9593},
   MRCLASS = {58G25 (52A40 53C20)},
  MRNUMBER = {683635},
MRREVIEWER = {Scott\ Wolpert},
       URL = {http://www.numdam.org/item?id=ASENS_1982_4_15_2_213_0},
}

@article{brooks1986spectral,
  title={The spectral geometry of a tower of coverings},
  author={Brooks, Robert},
  journal={Journal of Differential Geometry},
  volume={23},
  number={1},
  pages={97--107},
  year={1986},
  publisher={Lehigh University}
}

@incollection {brooks2006combinatorial,
    AUTHOR = {Brooks, Robert},
     TITLE = {Combinatorial problems in spectral geometry},
 BOOKTITLE = {Curvature and topology of {R}iemannian manifolds ({K}atata,1985)},
    SERIES = {Lecture Notes in Math.},
    VOLUME = {1201},
     PAGES = {14--32},
 PUBLISHER = {Springer, Berlin},
      YEAR = {1986},
      ISBN = {3-540-16770-6},
   MRCLASS = {58G25 (11F72)},
  MRNUMBER = {859574},
MRREVIEWER = {P.\ G.\ Zograf},
       DOI = {10.1007/BFb0075645},
       URL = {https://doi.org/10.1007/BFb0075645},
}

@book {chavel2001isoperimetric,
    AUTHOR = {Chavel, Isaac},
     TITLE = {Isoperimetric inequalities},
    SERIES = {Cambridge Tracts in Mathematics},
    VOLUME = {145},
      NOTE = {Differential geometric and analytic perspectives},
 PUBLISHER = {Cambridge University Press, Cambridge},
      YEAR = {2001},
     PAGES = {xii+268},
      ISBN = {0-521-80267-9},
   MRCLASS = {58J35 (35B05 49Q15 53C20 58-02)},
  MRNUMBER = {1849187},
MRREVIEWER = {Thierry\ Coulhon},
}

@article{mantuano2005discretization,
  title={Discretization of compact {R}iemannian manifolds applied to the spectrum of {L}aplacian},
  author={Mantuano, Tatiana},
  journal={Annals of Global Analysis and Geometry},
  volume={27},
  number={1},
  pages={33--46},
  year={2005},
  publisher={Springer}
}

@article {colbois2016steklov,
    AUTHOR = {Colbois, Bruno and Girouard, Alexandre and Raveendran, Binoy},
     TITLE = {The {S}teklov spectrum and coarse discretizations of manifolds
              with boundary},
   JOURNAL = {Pure Appl. Math. Q.},
  FJOURNAL = {Pure and Applied Mathematics Quarterly},
    VOLUME = {14},
      YEAR = {2018},
    NUMBER = {2},
     PAGES = {357--392},
      ISSN = {1558-8599,1558-8602},
   MRCLASS = {53C21 (58C40)},
  MRNUMBER = {4047402},
MRREVIEWER = {Changwei\ Xiong},
       DOI = {10.4310/pamq.2018.v14.n2.a3},
       URL = {https://doi.org/10.4310/pamq.2018.v14.n2.a3},
}

@article{buser1984bipartition,
  title={On the bipartition of graphs},
  author={Buser, Peter},
  journal={Discrete applied mathematics},
  volume={9},
  number={1},
  pages={105--109},
  year={1984},
  publisher={North-Holland}
}

@incollection {colbois2006spectral,
    AUTHOR = {Colbois, Bruno},
     TITLE = {The spectrum of the {L}aplacian: a geometric approach},
 BOOKTITLE = {Geometric and computational spectral theory},
    SERIES = {Contemp. Math.},
    VOLUME = {700},
     PAGES = {1--40},
 PUBLISHER = {Amer. Math. Soc., Providence, RI},
      YEAR = {2017},
      ISBN = {978-1-4704-2665-1},
   MRCLASS = {58J50 (35J10 35P15 35R30)},
  MRNUMBER = {3748520},
       DOI = {10.1090/conm/700/14181},
       URL = {https://doi.org/10.1090/conm/700/14181},
}

@article {yu2024monotonicity,
    AUTHOR = {Yu, Chengjie and Yu, Yingtao},
     TITLE = {Monotonicity of {S}teklov eigenvalues on graphs and
              applications},
   JOURNAL = {Calc. Var. Partial Differential Equations},
  FJOURNAL = {Calculus of Variations and Partial Differential Equations},
    VOLUME = {63},
      YEAR = {2024},
    NUMBER = {3},
     PAGES = {Paper No. 79, 22},
      ISSN = {0944-2669,1432-0835},
   MRCLASS = {05C50 (39A12)},
  MRNUMBER = {4721070},
MRREVIEWER = {Vladimir\ Mityushev},
       DOI = {10.1007/s00526-024-02683-y},
       URL = {https://doi.org/10.1007/s00526-024-02683-y},
}

@article{shi2022comparison,
  title={Comparison of {S}teklov eigenvalues and {L}aplacian eigenvalues on graphs},
  author={Shi, Yongjie and Yu, Chengjie},
  journal={Proceedings of the American Mathematical Society},
  volume={150},
  number={4},
  pages={1505--1517},
  year={2022}
}

@article{perrin2019lower,
  title={Lower bounds for the first eigenvalue of the {S}teklov problem on graphs},
  author={Perrin, H{\'e}l{\`e}ne},
  journal={Calculus of Variations and Partial Differential Equations},
  volume={58},
  number={2},
  pages={67},
  year={2019},
  publisher={Springer}
}

@article{perrin2021isoperimetric,
  title={Isoperimetric upper bound for the first eigenvalue of discrete {S}teklov problems},
  author={Perrin, H{\'e}l{\`e}ne},
  journal={The Journal of Geometric Analysis},
  volume={31},
  number={8},
  pages={8144--8155},
  year={2021},
  publisher={Springer}
}

@article{he2022upper,
  title={Upper bounds for the {S}teklov eigenvalues on trees},
  author={He, Zunwu and Hua, Bobo},
  journal={Calculus of Variations and Partial Differential Equations},
  volume={61},
  number={3},
  pages={101},
  year={2022},
  publisher={Springer}
}

@article{shi2025comparisons,
  title={Comparisons of {D}irichlet, {N}eumann and {L}aplacian eigenvalues on graphs and their applications},
  author={Shi, Yongjie and Yu, Chengjie},
  journal={Calculus of Variations and Partial Differential Equations},
  volume={64},
  number={9},
  pages={1--22},
  year={2025},
  publisher={Springer}
}

@article{he2021steklov,
  title={Steklov flows on trees and applications},
  author={He, Zunwu and Hua, Bobo},
  journal={arXiv preprint arXiv:2103.07696},
  year={2021}
}

@article{lin2024maximize,
  title={Maximize the {S}teklov eigenvalue of trees},
  author={Lin, Huiqiu and Zhao, Da},
  journal={Bulletin of the London Mathematical Society},
  year={2025},
  publisher={Wiley Online Library}
}

@article{shi2022lichnerowicz,
  title={A {L}ichnerowicz-type estimate for {S}teklov eigenvalues on graphs and its rigidity},
  author={Shi, Yongjie and Yu, Chengjie},
  journal={Calculus of Variations and Partial Differential Equations},
  volume={61},
  number={3},
  pages={98},
  year={2022},
  publisher={Springer}
}

@article {friedman1993some,
    AUTHOR = {Friedman, Joel},
     TITLE = {Some geometric aspects of graphs and their eigenfunctions},
   JOURNAL = {Duke Math. J.},
  FJOURNAL = {Duke Mathematical Journal},
    VOLUME = {69},
      YEAR = {1993},
    NUMBER = {3},
     PAGES = {487--525},
      ISSN = {0012-7094,1547-7398},
   MRCLASS = {05C50 (58G25)},
  MRNUMBER = {1208809},
MRREVIEWER = {J\'ozef\ Dodziuk},
       DOI = {10.1215/S0012-7094-93-06921-9},
       URL = {https://doi.org/10.1215/S0012-7094-93-06921-9},
}

@article{lin2025estimates,
  title={Estimates of the first {D}irichlet eigenvalue of graphs},
  author={Lin, Huiqiu and Liu, Lianping and You, Zhe and Zhao, Da},
  journal={arXiv preprint arXiv:2510.04557},
  year={2025}
}

@article {MR1041245,
    AUTHOR = {Grone, Robert and Merris, Russell and Sunder, V. S.},
     TITLE = {The {L}aplacian spectrum of a graph},
   JOURNAL = {SIAM J. Matrix Anal. Appl.},
  FJOURNAL = {SIAM Journal on Matrix Analysis and Applications},
    VOLUME = {11},
      YEAR = {1990},
    NUMBER = {2},
     PAGES = {218--238},
      ISSN = {0895-4798},
   MRCLASS = {05C50 (05C05 15A18 15A42)},
  MRNUMBER = {1041245},
MRREVIEWER = {Bojan\ Mohar},
       DOI = {10.1137/0611016},
       URL = {https://doi.org/10.1137/0611016},
}

\end{document}